\documentclass[sn-mathphys-num]{sn-jnl}

\usepackage{graphicx}%
\usepackage{multirow}%
\usepackage{amsmath,amssymb,amsfonts}%
\usepackage{amsthm}%
\usepackage{mathrsfs}%
\usepackage[title]{appendix}%
\usepackage{xcolor}%
\usepackage{textcomp}%
\usepackage{manyfoot}%
\usepackage{booktabs}%
\usepackage{algorithm}%
\usepackage{algorithmicx}%
\usepackage{algpseudocode}%
\usepackage{listings}%

\theoremstyle{thmstyleone}%
\newtheorem{theorem}{Theorem}
\newtheorem{proposition}[theorem]{Proposition}%

\theoremstyle{thmstyletwo}%
\newtheorem{example}{Example}%

\theoremstyle{thmstylethree}%
\newtheorem{definition}{Definition}%

\def\R{\mathbb{R}}
\def\K{\mathbb{K}}

\newcommand{\sgn}{\operatorname{sgn}}

\begin{document}

\title[Stationarity in nonsmooth optimization between geometrical motivation and topological relevance]{Stationarity in nonsmooth optimization between geometrical motivation and topological relevance}


\author[1]{\fnm{Vladimir} \sur{Shikhman}}\email{vladimir.shikhman@mathematik.tu-chemnitz.de}

\affil[1]{
\orgname{Chemnitz University of Technology}, \orgaddress{\street{Reichenhainer Str. 41}, \city{Chemnitz}, \postcode{09126}, \country{Germany}}}


\abstract{The goal of this paper is to compare alternative stationarity notions in structured nonsmooth optimization (SNO). Here, nonsmoothness is caused by complementarity, vanishing, orthogonality type, switching, or disjunctive constraints. On one side, we consider geometrically motivated notions of $\widehat N$-, $N$-, and $\overline{N}$-stationarity in terms of Fr\'echet, Mordukhovich, and Clarke normal cones to the feasible set, respectively. On the other side, we advocate the notion of topologically relevant T-stationarity, which adequately captures the global structure of SNO. Our main findings say that (a) $\widehat N$-stationary points include all local minimizers; (b) $N$-stationary points, which are not $\widehat N$-stationary, correspond to the singular saddle points of first order; (c) T-stationary points, which are not $N$-stationary, correspond to the regular saddle points of first order;  (d) $\overline{N}$-stationary points, which are not T-stationary, are irrelevant for optimization purposes. Overall, a hierarchy of stationarity notions for SNO is established.    

}




\maketitle

\section{Introduction}

Often nonsmooth structures enter into optimization with the description of the feasible set. Here, we consider mathematical programs with complementarity (MPCC), vanishing (MPVC), orthogonality type (MPOC), switching (MPSC), and disjunctive (MPDC) constraints. Our focus will be on the hierarchy of alternative stationarity notions for these subclasses of structured nonsmooth optimization (SNO). In nonlinear programming, where the feasible set is given by smooth equality and inequality constraints, the notion of Karush-Kuhn-Tucker points has meanwhile become standard.
In contrast to that, quite different stationarity notions were proposed in the framework of SNO in recent decades. All of them, while having certain justifications and being in use, can be said to currently compete throughout the literature. In particular, we consider geometrically motivated notions of $\widehat N$-, $N$-, and $\overline{N}$-stationarity. They are stated in terms of Fr\'echet, Mordukhovich, and Clarke normal cones to the feasible set, respectively, see \cite{rockafellar:1998, mordukhovich:2005}. In opposition to the latter, we advocate the topologically relevant notion of T-stationarity pioneered by Jongen, see  \cite{jongen:2000, shikhman:2012}. The concept of T-stationarity allows to adequately describe the global structure of SNO. For that, the topological changes of lower level sets while passing a T-stationary point are studied. First, we relate geometrically motivated and topologically relevant stationarity notion to each other by showing the following implications:
\[
    \widehat N\mbox{-stationarity} \,\,  \Longrightarrow \,\,N\mbox{-stationarity}  \,\,  \Longrightarrow \,\, \mbox{T-stationarity}
   \,\,  \Longrightarrow \,\, \overline{N}\mbox{-stationarity}.
\]
Concerning our main findings for SNO, we conclude that
\begin{itemize}
    \item $\widehat N$-stationary points include all local minimizers;
    \item $N$-, but not $\widehat N$-stationary points are singular saddle points of first order;
    \item T-, but not $N$-stationary points are regular saddle points of first order;
    \item $\overline{N}$-, but not T-stationary points are irrelevant for optimization purposes.
\end{itemize}
For these results to be stated, the novel notion of singular/regular saddle points of first order becomes crucial. Roughly speaking, a regular saddle point of first order admits two directions within the underlying nonsmooth structure, along which the objective function linearly decreases. In other words, it is characterized by having a non-vanishing biactive index. In turn, a singular saddle points of first order can be made regular by arbitrarily small $C^1$-perturbations of the SNO defining functions. We emphasize that the relations between geometrically motivated and topologically relevant stationarity notions are now fully clarified for SNO. It follows from our considerations that T-stationarity for nonsmooth optimization turns out to be the adequate generalization of the notion of Karush-Kuhn-Tucker points in the smooth case. For this conclusion, recall that saddle points are among the Karush-Kuhn-Tucker points, but they are exclusively of second order or higher. In case of nonlinear programming (NLP), the objective function has then a direction of quadratic decrease at such a saddle point. If facing nonsmooth structures, also saddle points of first order become relevant.  Actually, they cannot be  avoided in general if taking the perspective of global optimization. While $N$-stationarity adds singular saddle points of first order, only T-stationarity incorporates also their regular counterparts. 

The paper is organized as follows. In Section \ref{sec:sno}, the class of SNO is described. Sections \ref{sec:geom-stat} and \ref{sec:t-stat} are devoted to the introduction of geometrically motivated and topologically relevant stationarity notions for SNO, respectively. In Section \ref{sec:rel}, relations between these stationarity notions are presented.  

Our notation is standard. The cardinality of a finite set $A$ is denoted by $|A|$. The sign of a real number $a \in \R$ is denoted by $\sgn(a)$.
The $n$-dimensional
Euclidean space is denoted by $\mathbb{R}^n$. 
Given a twice continuously differentiable function $f:\mathbb{R}^n\rightarrow \mathbb{R}$, $\nabla f$ denotes its gradient, and $D^2f$ stands for its Hessian. 

\section{Structured nonsmooth optimization}
\label{sec:sno}

\begin{figure}[h]
    \centering
    \includegraphics[trim=600 570 700 85, scale=0.75]{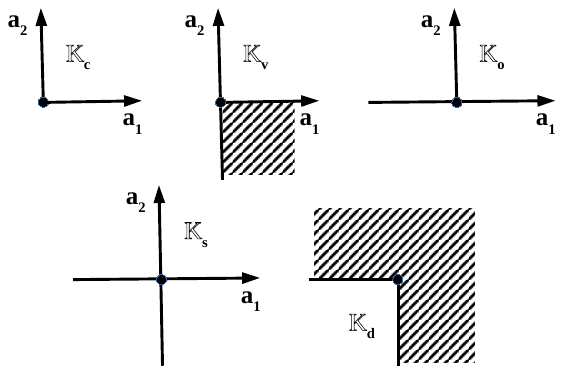}  
    \caption{Nonsmooth structures}
    \label{fig:cones}
\end{figure}

We consider the class of structured nonsmooth optimization problems:
\[
\mbox{SNO}: \quad
\min_{x} \,\, f(x)\quad \mbox{s.\,t.} \quad x \in M
\]
with
\[
    M=\left\{x \in\R^n \left\vert\;
       \left(F_{1,i}(x), F_{2,i}(x)\right)\in \K, i=1,\ldots,m 
     \right.\right\},
\]
were the defining functions $f, F_{1,i}, F_{2,i}:\R^n \rightarrow \R$, $i=1,\ldots,m$, are twice continuously differentiable, and the closed, but nonconvex cone $\K \subset \R^2$  represents nonsmooth structures to be studied, e.g. complementarity, vanishing, orthogonality type, switching, or disjunctive constraints. Accordingly, we restrict our attention to the following subclasses of SNO extensively discussed in the literature:
\begin{itemize}
    \item mathematical programs with complementarity constraints (MPCC), see e.g. \cite{ye:2005, jongen:2009}, with
\[
\K_\text{c}=\left\{ (a_1,a_2) \in \R^2\, \left\vert\, a_1\cdot a_2 = 0, a_1 \geq 0, a_2 \geq 0\right.\right\};
\]    
 \item mathematical programs with vanishing constraints (MPVC), see e.g. \cite{hoheisel:2007, dorsch:2012}, with
\[
\K_\text{v}=\left\{ (a_1,a_2) \in \R^2\, \left\vert\, a_1 \geq 0, a_1\cdot a_2 \leq 0, \right.\right\};
\] 
 \item mathematical programs with orthogonality type constraints (MPOC), see e.g. \cite{laemmel:2023}, with
\[
\K_\text{o}=\left\{ (a_1,a_2) \in \R^2\, \left\vert\, a_1\cdot a_2 = 0, a_2 \geq 0\right.\right\};
\]
 \item mathematical programs with switching constraints (MPSC), see e.g. \cite{mehlitz:2020, shikhman:2022}, with
\[
\K_\text{s}=\left\{ (a_1,a_2) \in \R^2\, \left\vert\, a_1\cdot a_2 = 0\right.\right\};
\]
 \item mathematical programs with disjunctive constraints (MPDC), see e.g. \cite{jongen:1997, mehlitz:2019}, with
\[
  \K_\text{d}=\left\{ (a_1,a_2) \in \R^2\, \left\vert\, a_1 \geq 0 \text{ or } a_2 \geq 0\right.\right\}.
\]
\end{itemize}
In this paper, $\K$ therefore stands as a substitute for one of the cones $\K_\text{c}, \K_\text{v}, \K_\text{o}, \K_\text{s}$, or $\K_\text{d}$, see also Figure \ref{fig:cones} for an illustration.

Now, we recall a suitable constraint qualification for the above SNO subclasses.

\begin{definition}[SNO-LICQ]
The SNO-tailored linear independence constraint qualification (SNO-LICQ) is said to be fulfilled at $\bar x \in M$ if the following vectors are linearly independent:
\[
\begin{array}{l}
\nabla F_{1,i}(\bar x) \mbox{ for } i \in \{1, \ldots, k\} \mbox{ with } F_{1,i}(\bar x) =0, \\
\nabla F_{2,i}(\bar x) \mbox{ for } i \in \{1, \ldots, k\} \mbox{ with } F_{2,i}(\bar x) =0.
\end{array}
\]
\end{definition}

It is well-known that SNO-LICQ is not too restrictive. In order to explain in which sense this is the case, we use the strong $C^2$-topology on the space $C^2(\R^n,\R)$ of twice differentiable functions. Denoted by $C^2_s$, it is sometimes referred to as the Whitney-topology, cf. \cite{jongen:2000, hirsch:1976}. The basis of 
$C^2_s$-topology is given by
$C^2$-perturbations of a twice differentiable function which are controlled by means of a continuous positive function. On the space of SNO defining functions $(f,F_{1,i},F_{2,i}, i =1,\ldots,m) \in C^2\left(\R^n,\R^{1+2m}\right)$ we consider the corresponding product $C^2_s$-topology.

\begin{proposition}[Genericity of SNO-LICQ]
\label{prop:licq-gen}
Let $\mathcal{F}\subset C^2\left(\R^n,\R^{1+2m}\right)$
denote the subset of SNO defining functions for which SNO-LICQ holds at all feasible points. Then, $\mathcal{F}$ is
$C_s^2$-open and -dense. 
\end{proposition}
Motivated by Proposition \ref{prop:licq-gen}, we always assume that SNO-LICQ holds at a feasible point $\bar x \in M$ under consideration. Since we are interested in studying nonsmooth structures, let us also assume that at $\bar x \in M$ all constraints are biactive, i.e. 
\[
  I_{00}(\bar x) = \{1, \ldots,m\},
\]
where the biactive index set is denoted by
\[
  I_{00}(\bar x) = \left\{ i \in \{1,\ldots,m\} \,\left\vert\, F_{1,i}(\bar x) =0, F_{2,i}(\bar x) =0\right.\right\}.
\]
Otherwise, the validity of SNO-LICQ at $\bar x$ allows us to equivalently rewrite the $i$-th constraint $\left(F_{1,i}(x), F_{2,i}(x)\right)\in \K$ as a smooth equality or inequality constraint, at least locally at $\bar x$. Then, we can treat the $i$-th constraint within the standard optimization theory for NLP, see e.g. \cite{jongen:2000}. Note that geometrically motivated and topologically relevant stationarity concepts coincide for NLP. Both reduce to the Karush-Kuhn-Tucker points. 
For simplicity of exposition, we thus exclude this case from our considerations without loss of generality.  

\section{Geometrically motivated stationarity}
\label{sec:geom-stat}

 In our geometrical study, we mainly follow the standard references \cite{mordukhovich:2005, rockafellar:1998} on the variational analysis. Let us start  by defining normal vectors to the SNO feasible set.
\begin{definition}[Normal cones]
Let $\bar x \in M$ be a SNO feasible point.
\begin{itemize}
        \item Fr\'echet normal cone at $\bar x$ to $M$ consists of vectors $v \in \widehat N_M(\bar x)$ iff
        \[
            \langle v,w\rangle \leq 0 \mbox{ for all } w \in T_M(\bar x),
        \]
        where the tangential cone consists of vectors $w \in T_M(\bar x)$ iff
        \[
            \mbox{there exist } x^k \underset{M}{\rightarrow} \bar x \mbox{ and } \tau^k \searrow 0 \mbox{ with } \frac{x^k - \bar x}{\tau^k} \rightarrow w.   
        \]
        \item Mordukhovich (limiting) normal cone at $\bar x$ to $M$ consists of vectors $v \in N_M(\bar x)$ iff 
        \[
            \mbox{ there exist } x^k \underset{M}{\rightarrow} \bar x \mbox{ and } v^k \in \widehat N_M(\bar x) \mbox{ with } v^k \rightarrow v.
        \]
        \item Clarke normal cone at $\bar x$ to $M$ is defined as the convex closure
        \[
             \overline{N}_M(\bar x)=\mbox{cl con } N(\bar x).
        \]
    \end{itemize}
\end{definition}

It is clear that by the very definition we have
\begin{equation}
    \label{eq:geom-rel}
    \widehat N_M(\bar x) \subset
  N_M(\bar x) \subset \overline{N}_M(\bar x).
\end{equation}

Normal vectors to the feasible set $M$ can be expressed in terms of the corresponding cones $\K$. By virtue of SNO-LICQ at $\bar x \in M$ and recalling that all nonsmooth constraints are assumed to be biactive, \cite[Theorem 6.14]{rockafellar:1998} provides the following formulae:
\[
\begin{array}{rcl}
     \widehat N_M(\bar x) &=& \displaystyle \left\{ \sum_{i \in I_{00}(\bar x)}\alpha_{1,i} \nabla F_{1,i}(\bar x)+\alpha_{2,i} \nabla F_{2,i}(\bar x)\, \left \vert \begin{array}{l}
         \left(\alpha_{1,i}, \alpha_{2,i} \right) \in \widehat N_{\mathbb{K}}(0), i \in I_{00}(\bar x)
   \end{array}\right.\right\},
 \\ \\
   N_M(\bar x) &=& \displaystyle \left\{ \sum_{i \in I_{00}(\bar x)} \alpha_{1,i} \nabla F_{1,i}(\bar x)+\alpha_{2,i} \nabla F_{2,i}(\bar x)\, \left \vert \begin{array}{l}
         \left(\alpha_{1,i}, \alpha_{2,i} \right) \in  N_{\mathbb{K}}(0), i \in I_{00}(\bar x)
   \end{array}\right.\right\},
\\ \\
   \overline{N}_M(\bar x) &=&\displaystyle \left\{ \sum_{i \in I_{00}(\bar x)}\alpha_{1,i} \nabla F_{1,i}(\bar x)+\alpha_{2,i} \nabla F_{2,i}(\bar x)\, \left \vert \begin{array}{l}
         \left(\alpha_{1,i}, \alpha_{2,i} \right) \in \overline{N}_{\mathbb{K}}(0), i \in I_{00}(\bar x)
   \end{array}\right.\right\}.
\end{array}
 \]
It is straightforward to explicitly compute Fr\'echet, Mordukhovich, and Clarke normal cones to $\K_\text{c}, \K_\text{v}, \K_\text{o}, \K_\text{s}$, and $\K_\text{d}$ at the origin:
\begin{itemize}
    \item MPCC:
    \[
\begin{array}{rcl}
      \widehat N_{\mathbb{K}_\text{c}}(0)&=&\displaystyle \left\{ (\alpha_1,\alpha_2) \in \R^2\, \left\vert\, \alpha_1 \leq 0, \alpha_2 \leq 0\right.\right\}, \\
      N_{\mathbb{K}_\text{c}}(0)&=&\displaystyle \left\{ (\alpha_1,\alpha_2) \in \R^2\, \left\vert\, \alpha_1 < 0, \alpha_2 < 0 \mbox{ or } \alpha_1 \cdot \alpha_2 =0\right.\right\}, \\  
      \overline{N}_{\mathbb{K}_\text{c}}(0)&=&\R^2;\\
\end{array}
\]
    \item MPVC:
    \[
\begin{array}{rcl}
      \widehat N_{\mathbb{K}_\text{v}}(0)&=&\displaystyle \left\{ (\alpha_1,\alpha_2) \in \R^2\, \left\vert\, \alpha_1 \leq 0, \alpha_2 = 0\right.\right\}, \\
      N_{\mathbb{K}_\text{v}}(0)&=&\displaystyle \left\{ (\alpha_1,\alpha_2) \in \R^2\, \left\vert\, \alpha_1 \cdot \alpha_2 = 0, \alpha_2 \geq 0 \right.\right\}, \\  
      \overline{N}_{\mathbb{K}_\text{v}}(0)&=&\displaystyle \left\{ (\alpha_1,\alpha_2) \in \R^2\, \left\vert\, \alpha_2 \geq 0 \right.\right\}; \\
\end{array}
\]
    \item MPOC:
    \[
\begin{array}{rcl}
      \widehat N_{\mathbb{K}_\text{o}}(0)&=&\displaystyle \left\{ (\alpha_1,\alpha_2) \in \R^2\, \left\vert\, \alpha_1 = 0, \alpha_2 \leq 0\right.\right\}, \\
      N_{\mathbb{K}_\text{o}}(0)&=&\displaystyle \left\{ (\alpha_1,\alpha_2) \in \R^2\, \left\vert\, \alpha_1 \cdot \alpha_2 = 0 \right.\right\}, \\  
      \overline{N}_{\mathbb{K}_\text{o}}(0)&=&\R^2; \\
\end{array}
\]
    \item MPSC:
    \[
\begin{array}{rcl}
      \widehat N_{\mathbb{K}_\text{s}}(0)&=&\displaystyle \left\{(0,0)\right\}, \\
      N_{\mathbb{K}_\text{s}}(0)&=&\displaystyle \left\{ (\alpha_1,\alpha_2) \in \R^2\, \left\vert\, \alpha_1 \cdot \alpha_2 = 0 \right.\right\}, \\  
      \overline{N}_{\mathbb{K}_\text{s}}(0)&=&\R^2; \\
\end{array}
\]
    \item MPDC:
    \[
\begin{array}{rcl}
      \widehat N_{\mathbb{K}_\text{d}}(0)&=&\displaystyle \left\{ (0,0) \right\}, \\
      N_{\mathbb{K}_\text{d}}(0)&=&\displaystyle \left\{ (\alpha_1,\alpha_2) \in \R^2\, \left\vert\, \alpha_1 \cdot \alpha_2 = 0, \alpha_1 \leq 0,\alpha_2 \leq 0 \right.\right\}, \\  
      \overline{N}_{\mathbb{K}_\text{d}}(0)&=&\displaystyle \left\{ (\alpha_1,\alpha_2) \in \R^2\, \left\vert\, \alpha_1 \leq 0, \alpha_2 \leq 0 \right.\right\}. \\
\end{array}
\]
\end{itemize}

Now, we are ready to define geometrically motivated stationarity notions.
\begin{definition}[Geometrically motivated stationarity]
  A SNO feasible point $\bar x \in M$ is called 
  \begin{itemize}
      \item $\widehat N$-stationary if
    $0 \in \nabla f(\bar x) + \widehat N_M(\bar x)$,
   \item $N$-stationary if
   $0 \in \nabla f(\bar x) + N_M(\bar x)$,
   \item $\overline{N}$-stationary if
   $0 \in \nabla f(\bar x) + \overline{N}_M(\bar x)$. 
  \end{itemize}
\end{definition}

Due to (\ref{eq:geom-rel}), we can easily relate the geometrically motivated stationarity notions to each other as follows:
\begin{equation}
    \label{eq:impl-stat}
    \widehat N\mbox{-stationarity} \,\,  \Longrightarrow \,\,N\mbox{-stationarity}
   \,\,  \Longrightarrow \,\, \overline{N}\mbox{-stationarity}.
\end{equation}  
In presence of SNO-LICQ, geometrically motivated stationarity turns out to be necessary for optimality, see e.g. \cite[Theorem 6.12]{rockafellar:1998}.

\begin{proposition}[Necessary optimality condition]
    \label{prop:noc}
    Let $\bar x \in M$ be a local minimizer of SNO satisfying SNO-LICQ. Then, $\bar x$ is $\widehat N$-stationary and, thus, $N$- and $\overline{N}$-stationary as well.
\end{proposition}

\begin{figure}[h]
    \centering
    \includegraphics[trim=670 270 700 145, scale=0.75]{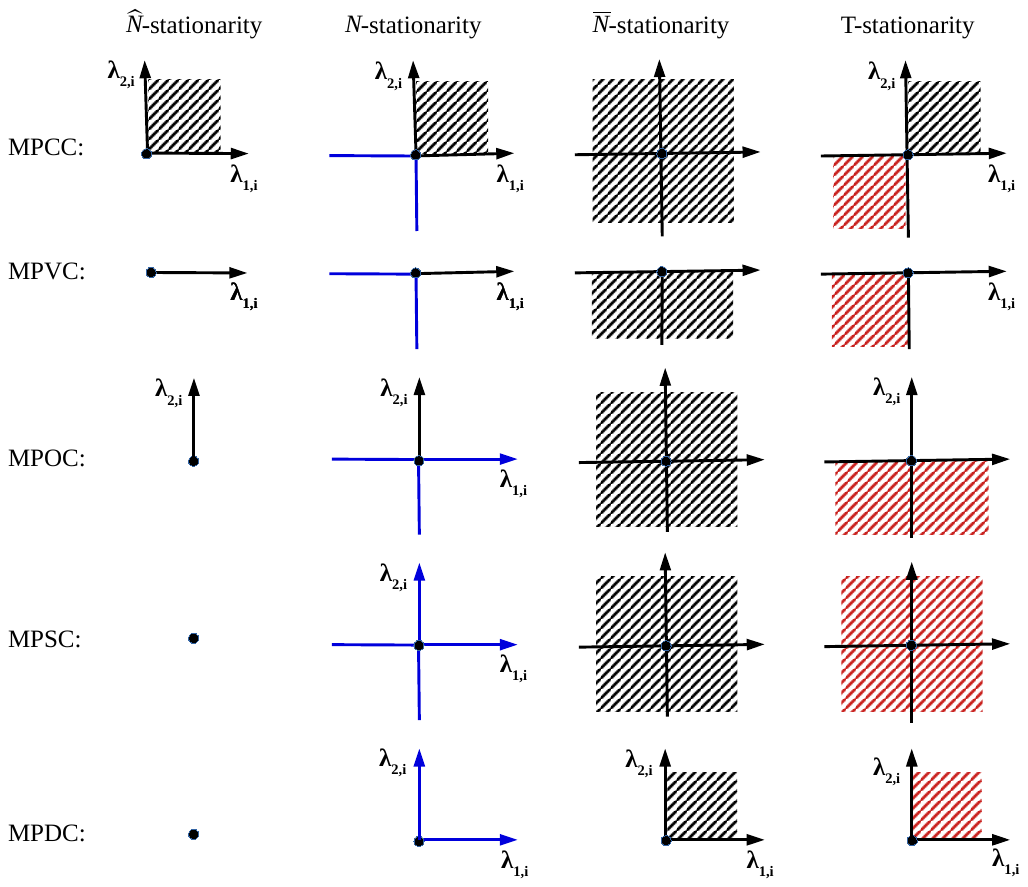}  
    \caption{Biactive multipliers $\lambda_{1,i}$ and $\lambda_{2,i}$ for $\widehat N$-, $N$-, $\overline{N}$-, and T-stationarity}
    \label{fig:stat1}
\end{figure}

It is convenient to rewrite geometrically motivated stationarity notions in the form
\begin{equation}
    \label{eq:stat-grad}
    \nabla f(\bar x) = \sum_{i \in I_{00}(\bar x)} \lambda_{1,i} \nabla F_{1,i}(\bar x) + \lambda_{2,i} \nabla F_{2,i}(\bar x),
\end{equation}
where the multipliers $\lambda_{1,i}$, $\lambda_{2,i}$, $i\in I_{00}(\bar x)$, are restricted correspondingly. By setting $\lambda_{1,i}=-\alpha_{1,i}$ and $\lambda_{2,i}=- \alpha_{2,i}$
in the formulae for the normal cones above, we obtain for the pairs of biactive multipliers  in (\ref{eq:stat-grad}) the following cases:
\begin{itemize}
    \item MPCC:
    \[
\begin{array}{rcl}
     \widehat N\mbox{-}stationarity & \Leftrightarrow& \lambda_{1,i}\geq 0, \lambda_{2,i} \geq 0, i\in I_{00}(\bar x), \\ 
    N\mbox{-}stationarity & \Leftrightarrow& \lambda_{1,i} > 0, \lambda_{2,i} > 0 \mbox{ or } \lambda_{1,i} \cdot \lambda_{2,i}=0, i\in I_{00}(\bar x), \\
     \overline{N}\mbox{-}stationarity & \Leftrightarrow& \mbox{no restrictions on }\lambda_{1,i}, \lambda_{2,i}, i\in I_{00}(\bar x); \\
\end{array}
\]
    \item MPVC:
    \[
\begin{array}{rcl}
     \widehat N\mbox{-}stationarity & \Leftrightarrow& \lambda_{1,i}\geq 0, \lambda_{2,i} = 0, i\in I_{00}(\bar x), \\ 
     N\mbox{-}stationarity & \Leftrightarrow& \lambda_{1,i} \cdot \lambda_{2,i} =0, \lambda_{2,i} \leq 0, i\in I_{00}(\bar x), \\
     \overline{N}\mbox{-}stationarity & \Leftrightarrow&    \lambda_{2,i} \leq 0, i\in I_{00}(\bar x); \\
\end{array}
\]
    \item MPOC:
    \[
\begin{array}{rcl}
     \widehat N\mbox{-}stationarity & \Leftrightarrow& \lambda_{1,i}= 0, \lambda_{2,i} \geq 0, i\in I_{00}(\bar x), \\ 
    N\mbox{-}stationarity & \Leftrightarrow&  \lambda_{1,i} \cdot \lambda_{2,i}=0, i\in I_{00}(\bar x), \\
     \overline{N}\mbox{-}stationarity & \Leftrightarrow& \mbox{no restrictions on }\lambda_{1,i}, \lambda_{2,i}, i\in I_{00}(\bar x); \\
\end{array}
\]
    \item MPSC:
    \[
\begin{array}{rcl}
     \widehat N\mbox{-}stationarity & \Leftrightarrow& \lambda_{1,i}= 0, \lambda_{2,i}= 0, i\in I_{00}(\bar x), \\ 
    N\mbox{-}stationarity & \Leftrightarrow&  \lambda_{1,i} \cdot \lambda_{2,i}=0, i\in I_{00}(\bar x), \\
     \overline{N}\mbox{-}stationarity & \Leftrightarrow& \mbox{no restrictions on }\lambda_{1,i}, \lambda_{2,i}, i\in I_{00}(\bar x); \\
\end{array}
\]
    \item MPVC:
    \[
\begin{array}{rcl}
     \widehat N\mbox{-}stationarity & \Leftrightarrow& \lambda_{1,i}= 0, \lambda_{2,i} = 0, i\in I_{00}(\bar x), \\ 
     N\mbox{-}stationarity & \Leftrightarrow& \lambda_{1,i} \cdot \lambda_{2,i} =0, \lambda_{1,i} \geq 0, \lambda_{2,i} \geq 0, i\in I_{00}(\bar x), \\
     \overline{N}\mbox{-}stationarity & \Leftrightarrow&    \lambda_{1,i} \geq 0, \lambda_{2,i} \geq 0, i\in I_{00}(\bar x). \\
\end{array}
\]
\end{itemize}
In Figure \ref{fig:stat1}, the notions of $\widehat N$-, $N$-, and $\overline{N}$-stationarity   are illustrated for the subclasses MPCC, MPVC, MPOC, MPSC, and MPDC.

Let us briefly comment on other related stationarity notions sometimes used for SNO in the literature:

\begin{itemize}
    \item S-stationarity stands for "strong" and is derived by putting SNO into the framework of NLP. This can be done when the SNO feasible set is given by smooth equality and inequality constraints. Then, the notion of Karush-Kuhn-Tucker points is applied even if the standard linear independence constraint qualification (LICQ) is normally violated. To overcome this obstacle, weaker constraint qualification  of e.g. Abadie or Guignard type are shown to hold under SNO-LICQ. Let us illustrate this procedure for MPCC, see e.g. \cite{ye:2005}, whose feasible set can be written as
    \[
    M_\text{c}=\left\{x \in\R^n \left\vert\;
       F_{1,i}(x) \cdot F_{2,i}(x)=0, F_{1,i}(x) \geq 0,  F_{2,i}(x) \geq 0, i=1,\ldots,m 
     \right.\right\}.
    \]
   For a corresponding Karush-Kuhn-Tucker point $\bar x \in M_\text{c}$ there exist multipliers $\lambda_i \geq 0$, $\lambda_{1,i} \geq 0$, $\lambda_{2,i} \geq 0$, $i \in I_{00}(\bar x)$, such that it holds:
      \begin{equation}
     \label{eq:grad-s}
     \nabla f(\bar x) = \sum_{i\in I_{00}(\bar x)} \lambda_{i} \nabla \left(F_{1,i} \cdot F_{2,i} \right)(\bar x)+ \lambda_{1,i} \nabla F_{1,i}(\bar x) + \lambda_{2,i} \nabla F_{2,i}(\bar x).
\end{equation}
    By recalling that all constraints are assumed to be biactive, we obtain:
    \[
    \nabla \left(F_{1,i} \cdot F_{2,i}\right)(\bar x) = \nabla F_{1,i}(\bar x) \cdot F_{2,i}(\bar x) + F_{1,i}(\bar x) \cdot \nabla F_{2,i}(\bar x) =0.
    \]
    Hence, (\ref{eq:grad-s}) becomes
    \[
     \nabla f(\bar x) = \sum_{i\in I_{00}(\bar x)} \lambda_{1,i} \nabla F_{1,i}(\bar x) + \lambda_{2,i} \nabla F_{2,i}(\bar x),
    \]
    where
    \[
      \lambda_{1,i} \geq 0, \lambda_{2,i} \geq 0, i \in I_{00}(\bar x).
    \]
    Since the biactive multipliers are nonnegative, we see that S-stationarity coincides with $\widehat N$-stationarity for MPCC. The same is true for MPVC, MPSC, and MPOC as the analogous derivations show. Note that the MPDC feasible set cannot be written in terms of smooth equality and inequality constraints. Nevertheless, it has been proposed in \cite{mehlitz:2019} to define S-stationarity for MPDC also via $\widehat N$-stationarity. Overall, we conclude that the notions of $S$-and $\widehat N$-stationarity are identical for SNO.   
    \item M-stationarity coincides with $N$-stationarity and resembles the fact that the normal cone due to Mordukhovich is considered.   
    \item W-stationarity is defined just by requiring (\ref{eq:stat-grad}) with no restrictions on biactive multipliers $\lambda_{1,i}$ and $\lambda_{2,i}$. This loose stationarity notion obviously does not in general take the particular type of nonsmoothness into account.
\end{itemize}

\section{Topologically relevant stationarity}
\label{sec:t-stat}

In our topological study, we mainly follow the references \cite{jongen:2000, shikhman:2012} and papers \cite{shikhman:2022, laemmel:2023, jongen:1997}, where the ideas of Morse theory were introduced into the optimization context.  
By doing so, the notion of T-stationarity plays a crucial role. It helps to adequately describe the topological properties of SNO lower level sets, where $a\in \R$ is varying:
\[
M_a=\left\{x\in M \left\vert f(x)\le a\right.\right\}.
\]
Eventually, we also consider intermediate sets, which will be denoted for $a < b$ by
\[
M^b_a=\left\{x\in M \left\vert a\le f(x) \le b\right. \right\}.
\]

\begin{definition}[T-stationarity]
A SNO feasible point $\bar x \in M$ is called T-stationary if 
\begin{equation}
     \label{eq:grad}
     \nabla f(\bar x) = \sum_{i\in I_{00}(\bar x)} \lambda_{1,i} \nabla F_{1,i}(\bar x) + \lambda_{2,i} \nabla F_{2,i}(\bar x)
\end{equation}
holds with biactive multipliers  $\lambda_{1,i}$, $\lambda_{2,i}$, $i\in I_{00}(\bar x)$, which -- depending on the type of nonsmoothness --  are restricted as follows:
\begin{itemize}
    \item MPCC:  $\lambda_{1,i} \cdot \lambda_{2,i} \geq 0, i\in I_{00}(\bar x)$,
    \item MPVC:  $\lambda_{1,i} \cdot \lambda_{2,i} \geq 0, \lambda_{2,i} \leq 0$, $i\in I_{00}(\bar x)$,
    \item MPOC:  $\lambda_{1,i}=0 \mbox{ or }\lambda_{2,i} \leq 0, i\in I_{00}(\bar x)$,
    \item MPSC:  no restrictions on $\lambda_{1,i}, \lambda_{2,i}, i\in I_{00}(\bar x)$,
    \item MPDC:  $\lambda_{1,i} \geq 0, \lambda_{2,i} \geq 0, i\in I_{00}(\bar x)$.    
\end{itemize}
\end{definition}

Again, Figure \ref{fig:stat1} provides a graphical illustration of T-stationarity for the subclasses
MPCC, MPVC, MPOC, MPSC, and MPDC.
In presence of SNO-LICQ, geometrically motivated stationarity turns out to be necessary for optimality as well.
\begin{proposition}[Necessary optimality condition]
\label{prop:necess}
Let $\bar x \in M$ be a local minimizer of SNO satisfying SNO-LICQ, then $\bar x$ is a T-stationary point.
\end{proposition}

Outside the set of T-stationary points, the topology of SNO lower level sets remains unchanged. This is referred to as deformation result within the scope of Morse theory, cf. \cite{jongen:2000}. 

\begin{theorem}[Deformation]
\label{thm:def}
Let $M^b_a$ be compact and SNO-LICQ be fulfilled at all points $x \in M_a^b$. Then, if the intermediate set $M_a^b$ contains no T-stationary points for SNO, then the lower level sets $M_a$ and $M_b$ are homeomorphic.
\end{theorem}

In order to describe the topological changes of  SNO lower level sets if crossing T-stationary points, we need to introduce the notion of their nondegeneracy. 
For that, given a T-stationary point $\bar x \in M$ with the multipliers $\lambda_{1,i}, \lambda_{2,i}$, $i\in I_{00}(\bar x)$, we consider the Lagrange function
\[
        \begin{array}{rcl}
        L(x)&=&\displaystyle   
        f(x)-
\sum_{i\in I_{00}(\bar x)} \lambda_{1,i}  F_{1,i}(x) + \lambda_{2,i}  F_{2,i}(x). \end{array}
\]
 We further use the corresponding tangent space
\[
\mathcal{T}_{\bar x}=\left\{\xi \in \R^n \,\left\vert \,
\begin{array}{l}
DF_{1,i}\left(\bar x\right)\xi=0, DF_{2,i}\left(\bar x\right)\xi=0, i \in I_{00}(\bar x)
\end{array}
\right.\right\}.
\]

\begin{definition}[Nondegeneracy]
\label{def:nondeg}
A T-stationary point $\bar x \in M$ of SNO with multipliers $\lambda_{1,i}, \lambda_{2,i}$, $i\in I_{00}(\bar x)$ is called nondegenerate if the following conditions are satisfied:
\begin{itemize}
    \item []ND1: SNO-LICQ holds at $\bar x$,
    \item []ND2: the biactive multipliers do not vanish, i.\,e. $\lambda_{1,i}\ne0,\lambda_{2,i}\ne 0$, $i\in I_{00}(\bar x)$,
    \item []ND3: the matrix $D^2 L(\bar x)\restriction_{T_{\bar x}}$ is nonsingular.
\end{itemize}
\end{definition}


It turns out that nondegeneracy holds at all T-stationary points of a generic SNO.

\begin{proposition}
[Nondegeneracy is generic]
\label{prop:generic}
Let $\mathcal{G}\subset C^2\left(\R^n,\R^{1+2k}\right)$ be the subset of SNO defining functions for which each 
T-stationary point is nondegenerate. Then, $\mathcal{G}$ is $C^2_s$-open and -dense.
\end{proposition}

With a T-stationary point it is further  convenient to associate its T-index. 

\begin{definition}[T-index]
Let $\bar x \in M$ be a T-stationary point of SNO satisfying SNO-LICQ and, thus, with unique multipliers $\lambda_{1,i}$, $\lambda_{2,i}$, $i\in I_{00}(\bar x)$. The  number of negative eigenvalues of the matrix $D^2 L(\bar x)\restriction_{\mathcal{T}_{\bar x}}$ is called its quadratic index $QI$. The biactive index $BI$ of $\bar x$ depends on the nonsmooth structure as follows:
\begin{itemize}
    \item MPCC: it is the number of negative biactive multipliers, i.e.
    \[
      BI = \left| \left\{i \in I_{00}(\bar x) \, \left\vert\, \lambda_{1,i} <0, \lambda_{2,i} <0\right.\right\}\right|;
    \]
    \item MPVC, MPOC, MPSC, or MPDC: it is the number of non-vanishing biactive multipliers, i.e.
    \[
      BI = \left|\left\{i \in I_{00}(\bar x) \, \left\vert\, \lambda_{1,i} \not = 0, \lambda_{2,i} \not =0\right.\right\}\right|.
    \]
\end{itemize}
We define the T-index as the sum of both, i.\,e. $TI=QI+BI$.
\end{definition}



Sufficient optimality condition for generic SNO can be derived in terms of T-index.

\begin{proposition}[Sufficient optimality condition]
\label{prop:zero-index}
A nondegenerate T-stationary point of SNO is a local minimizer if and only if
its T-index vanishes.
\end{proposition}

Proposition \ref{prop:zero-index} says that at a nondegenerate minimizer $\bar x \in M$ of SNO we have both quadratic and biactive indices vanishing. Note that $QI=0$ means that the matrix $D^2 L(\bar x)\restriction_{\mathcal{T}_{\bar x}}$ is positive definite. This is the usual second-order sufficient condition (SOSC). $BI=0$ means that in case of MPCC all biactive multipliers are positive, i.e.
\[
    \lambda_{1,i} > 0, \lambda_{2,i}>0, i \in I_{00}(\bar x),
\]
whereas in case of MPVC, MPOC, MPSC, or MPDC the biactive index set is empty, i.e.
\[
    I_{00}(\bar x) = \emptyset.
\]
How is it compatible with our assumption that all constraints are biactive at a point of interest? For MPVC, MPOC, MPSC, or MPDC, we conclude that a nondegenerate minimizer cannot have biactive constraints. The nondegenerate minimizers for those subclasses are just usual Karush-Kuhn-Tucker points then.   

In general, we emphasize that nondegenerate T-stationary points with non-vanishing T-index are saddle points for SNO. The local structure of SNO in the vicinity of nondegenerate saddle points is fully described by their quadratic and biactive indices, at least up to the smooth change of coordinates.  

\begin{proposition}[Morse Lemma]
\label{prop:morse}
Suppose that $\bar x \in M$ is a nondegenerate T-stationary point of SNO with multipliers $\lambda_{1,i}, \lambda_{2,i}$, $i\in I_{00}(\bar x)$. Let its quadratic index be $QI$ and its biactive index be $BI$. Then, there exist neighborhoods $U$ and $V$ of $\bar x$ and $0$, respectively, and a local $C^1$-coordinate system $\Psi: U \rightarrow V$ of $\R^n$ around $\bar x$ such that
\begin{equation}
\label{eq:normal}
    f\circ \Psi^{-1}(y)= f(\bar x) +
    \sum\limits_{i \in I_{00}(\bar x)}\left(\sgn(\lambda_{1,i}) \cdot y_{1,i}+\sgn(\lambda_{2,i}) \cdot y_{2,i}\right) + \sum\limits_{j=2\left|I_{00}(\bar x)\right|+1}^n\pm y_{j}^2,  
\end{equation}
where $y \in \mathbb{K}^{\left|I_{00}(\bar x)\right|} \times \R^{n-2\left|I_{00}(\bar x)\right|}$. Moreover, there are exactly $QI$ negative squares in (\ref{eq:normal}). In case of MPCC, the number of negative linear terms in (\ref{eq:normal}) is $BI$; in case of MPVC, MPOC, MPSC, or MPDC, the number of linear terms in (\ref{eq:normal}) is $BI$.
\end{proposition}

Now, we are ready to generally describe how the topology of the SNO lower level sets changes if passing a T-stationary level. 

\begin{theorem}[Cell-attachment]
\label{thm:cell-a}
Let $M_a^b$ be compact and suppose that it contains exactly one T-stationary point $\bar x$. Furthermore, let $\bar x$ be nondegenerate with T-index equal to $t$. 
If $a<f(\bar x) <b$,
then $M_b$ is homotopy-equivalent to $M_a$ with a $t$-cell attached along its boundary.
\end{theorem}

\begin{figure}[h]
    \centering
    \includegraphics[trim=580 470 700 5, scale=0.6]{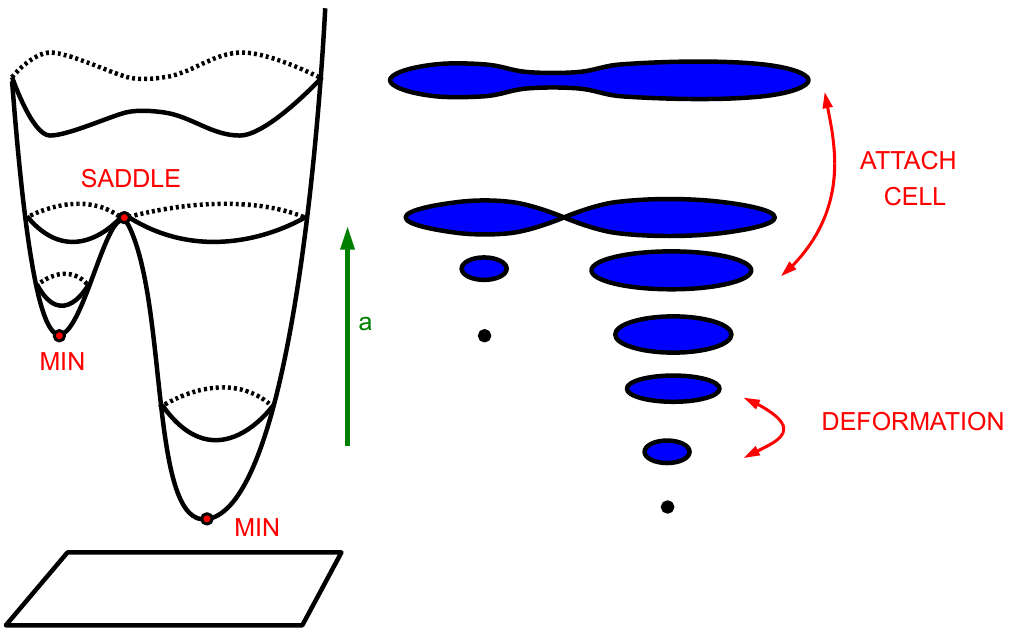}  
    \caption{Deformation and Cell-attachment}
    \label{fig:morse}
\end{figure}

A global interpretation of our results is typical for the Morse theory, cf. \cite{jongen:2000}. Suppose that the SNO feasible
set is compact and connected, that 
SNO-LICQ holds at all feasible points, and that all T-stationary points are nondegenerate with pairwise
different functional values. Then, a connected component of the lower level set is created as soon as we pass a level corresponding to a local minimizer. Different components can
only be connected if a $1$-cell is attached. This implies that there are at least $(k-1)$ T-stationary
points with T-index equal to one. Here, $k$ denotes the number of local minimizers of SNO. In the context of global optimization the latter result usually is referred to as mountain pass. For NLP, it is well known that Karush-Kuhn-Tucker points with quadratic index equal to one naturally appear along with local minimizers, see Figure \ref{fig:morse}. For SNO, however, not only T-stationary points with quadratic index equal to one, but also with biactive index equal to one may become relevant. 

At the end of this section we additionally relate T- to the so called C-stationary points. C-stationarity refers to Clarke's necessary optimality conditions, when the feasible set can be expressed by means of nonsmooth equality and inequality constraints, see \cite{clarke:1983}. In our setting, this is the case for MPCC with
    \[
    M_\text{c}=\left\{x \in\R^n \left\vert\;
      \min \left\{ F_{1,i}(x), F_{2,i}(x)\right\}=0, i=1,\ldots,m 
     \right.\right\}
    \]
and for MPDC with    
\[
    M_\text{d}=\left\{x \in\R^n \left\vert\;
      \max \left\{ F_{1,i}(x), F_{2,i}(x)\right\}\geq 0, i=1,\ldots,m 
     \right.\right\}.
    \]
Since the derivations of C-stationarity are similar in both cases, let us focus first on MPCC. 
    For a corresponding C-stationary point $\bar x \in M_\text{c}$ there exist multipliers $\lambda_i$, $i \in I_{00}(\bar x)$, such that it holds:
      \begin{equation}
     \label{eq:grad-c}
     \nabla f(\bar x) = \sum_{i\in I_{00}(\bar x)} \lambda_{i} \overline{\partial} \min \left\{ F_{1,i}, F_{2,i}\right\}(\bar x).
\end{equation}
    Here, by $\overline{\partial}g$ we denote the subdifferential due to Clarke of a locally Lipschitz continuous function $g:\R^n \rightarrow \R$. More precisely, it is defined as the convex hull
    \[
      \overline{\partial}g(\bar x)=\text{conv}\left\{ \lim \nabla g(x^k)\left\vert\; x^k \rightarrow \bar x, x^k \not\in \Omega_g\right.\,\right\},
    \]
    where $\Omega_g \subset \R^n$ denotes the set of points, at which $g$  fails to be differentiable.
    Due to \cite[Proposition 2.3.12]{clarke:1983}, the Clarke's subdifferential of the maximum of differentiable functions can be computed as
\[
\overline{\partial} \max \left\{ F_{1,i}, F_{2,i}\right\}(\bar x) = \mbox{conv} \left\{\nabla F_{1,i}(\bar x), \nabla F_{2,i}(\bar x)\right\}.
\]    
Recall for the latter that all constraints are assumed to be biactive at $\bar x$. Further, we use the property $\overline{\partial}(- \cdot)=-\overline{\partial}(\cdot)$ of the Clarke's subdifferential, see \cite[Proposition 2.3.1]{clarke:1983}, to derive
\[
   \overline{\partial} \min \left\{ F_{1,i}, F_{2,i}\right\}(\bar x) = 
   \mbox{conv} \left\{\nabla F_{1,i}(\bar x), \nabla F_{2,i}(\bar x)\right\}.
\]
Substituting into (\ref{eq:grad-c}), we obtain $\beta_{i} \in [0,1]$, $i \in I_{00}(\bar x)$, with
\[
       \nabla f(\bar x) = \sum_{i\in I_{00}(\bar x)} \lambda_{i} \left(\beta_i\nabla F_{1,i}(\bar x) + (1-\beta_i) \nabla F_{2,i}(\bar x)\right).
\]
By setting $\lambda_{1,i}=\lambda_i \cdot \beta_i$, $\lambda_{1,i}=\lambda_i \cdot (1-\beta_i)$, $i \in I_{00}(\bar x)$, we equivalently have
\begin{equation}
    \label{eq:stat-c1}
    \nabla f(\bar x) = \sum_{i\in I_{00}(\bar x)} \lambda_{1,i} \nabla F_{1,i}(\bar x) + \lambda_{2,i} \nabla F_{2,i}(\bar x),
\end{equation}
where the biactive multipliers are of the same sign, i.e.
\[
  \lambda_{1,i} \cdot \lambda_{2,i} \geq 0, i \in I_{00}(\bar x).
\]
We see that this is exactly the definition of T-stationarity for MPCC. 
Analogously for MPDC, C-stationarity amounts to 
     \begin{equation}
     \label{eq:grad-c1}
     \nabla f(\bar x) = \sum_{i\in I_{00}(\bar x)} \lambda_{i} \overline{\partial} \max \left\{ F_{1,i}, F_{2,i}\right\}(\bar x).
\end{equation}
with $\lambda_i \geq 0$, $i \in I_{00}(\bar x)$. The latter sign condition mimics the definition of Karush-Kuhn-Tucker points in the smooth setting, see \cite[Theorem 6.1.1]{clarke:1983}. After renaming the multipliers as above, we get (\ref{eq:stat-c1}) with \[
  \lambda_{1,i} \geq 0, \lambda_{2,i} \geq 0, i \in I_{00}(\bar x).
\]
From the nonnegativity of biactive multipliers we see that Clarke's necessary optimality conditions provide T-stationary points also for MPDC. We note that the concept of C-stationarity is not applicable for MPVC, MPOC, and MPSC. In fact, the corresponding feasible sets $M_\text{v}$, $M_\text{o}$, and $M_\text{s}$ cannot be naturally represented by nonsmooth equality and inequality constraints. At the origin, they are namely not homeomorhic to a Lipschitz manifold with boundary of a certain dimension. Therefore, the derivation of T-stationarity for MPVC, MPOC, and MPSC, but also for other potential SNO subclasses remains a challenging task. Although it coincides with C-stationarity for MPCC and MPDC, we doubt if any nonsmooth calculus can be used for deriving T-stationarity in general.

\section{Relations between stationarity notions}
\label{sec:rel}

We start by incorporating the topologically relevant stationarity into the scheme (\ref{eq:impl-stat}) of implications between geometrically motivated stationarity notions.
This can be easily done by comparing the signs of the corresponding biactive multipliers, cf. Figure \ref{fig:stat1}.
It turns out that T-stationarity is weaker than $N$-, but stronger than $\overline{N}$-stationarity.  

\begin{proposition}[General relations]
    Under SNO-LICQ, we have the following relations between geometrically motivated and topologically relevant stationarity notions:
    \begin{equation}
    \label{eq:impl-stat1}
    \widehat N\mbox{-stationarity} \,\,  \Longrightarrow \,\,N\mbox{-stationarity}  \,\,  \Longrightarrow \,\, \mbox{T-stationarity}
   \,\,  \Longrightarrow \,\, \overline{N}\mbox{-stationarity}.
\end{equation}  
   Moreover, none of the converse implications is true in general.
\end{proposition}

Now, we turn our attention to the differences between $\widehat N$-, $N$-, $\overline{N}$-, and T-stationarity. The question is what kind of stationary points should be added to a stronger notion in order to get a weaker one? 
For that, we first introduce the notion of saddle points of first order for SNO. It is motivated by the role the biactive index and the nondegeneracy condition ND2 play for T-stationarity within the Morse theory, see Section \ref{sec:t-stat}. We just relax ND2 by allowing at a saddle point of first order one of its biactive multipliers to vanish, but not both of them. Concerning the restrictions for biactive multipliers of saddle points of first order, we mimic the condition $BI \not =0$ by allowing for mentioned degeneracies.   

\begin{definition}[Saddle points of first order]
Let T-stationary point $\bar x \in M$ is called saddle point of first order if for some $i \in I_{00}(\bar x)$ the corresponding biactive multipliers -- depending on the type of
nonsmoothness -- are restricted as follows::
\begin{itemize}
    \item MPCC:  $\lambda_{1,i} \leq 0, \lambda_{2,i} \leq 0$ and $\left(\lambda_{1,i},\lambda_{2,i}\right) \not = (0,0)$,
    \item MPVC:  $\lambda_{1,i} \leq 0, \lambda_{2,i} \leq 0$ and $\left(\lambda_{1,i},\lambda_{2,i}\right) \not = (0,0)$,
    \item MPOC:  $\lambda_{2,i} \leq 0$ and $\left(\lambda_{1,i},\lambda_{2,i}\right) \not = (0,0)$, 
    \item MPSC:  $\left(\lambda_{1,i},\lambda_{2,i}\right) \not = (0,0)$, 
    \item MPDC:  $\lambda_{1,i} \geq 0, \lambda_{2,i} \geq 0$ and $\left(\lambda_{1,i},\lambda_{2,i}\right) \not = (0,0)$.    
\end{itemize}   
\end{definition}

\begin{figure}[h]
    \centering
    \includegraphics[trim=730 270 700 205, scale=0.75]{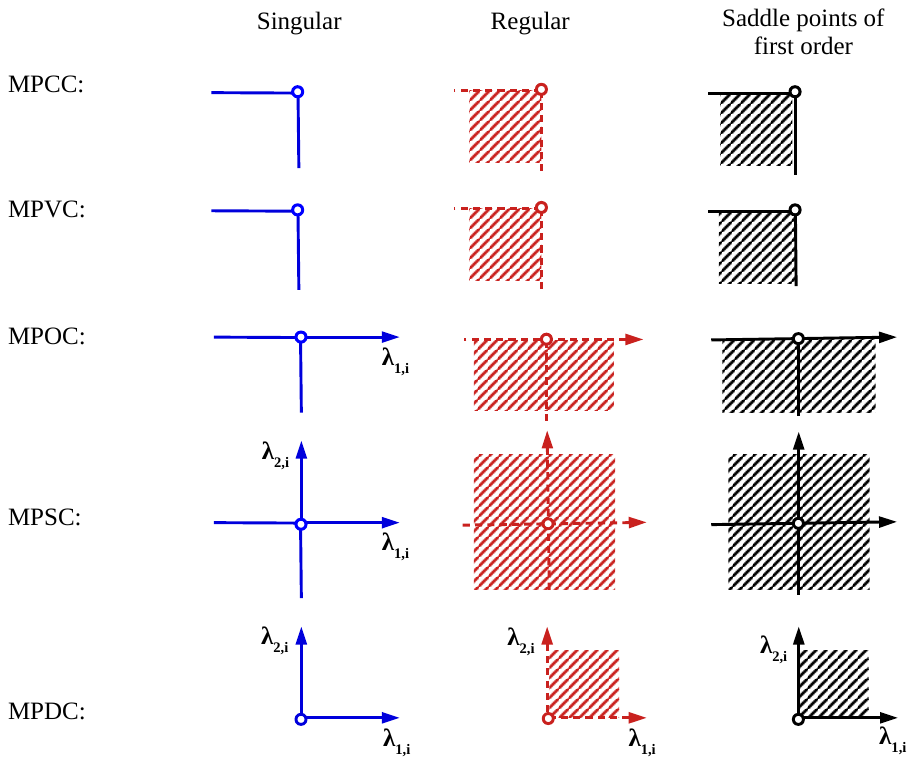}  
    \caption{Biactive multipliers $\lambda_{1,i}$ and $\lambda_{2,i}$ for singular (in blue) and regular (in red) saddle points of first order (in black)}
    \label{fig:spfo}
\end{figure}

The difference between T- and $\widehat N$-stationarity can be captured by means of the saddle points of first order, cf. Figures \ref{fig:stat1} and \ref{fig:spfo}.

\begin{proposition}[$\widehat N$- vs. T-stationarity]
\label{prop:t-hat-n}
  A T-stationary point $\bar x\in M$ satisfying SNO-LICQ is not $\widehat N$-stationary for SNO if and only if it is a saddle point of first order.   
\end{proposition}

Note that $\widehat N$-stationarity is a necessary optimality condition under SNO-LICQ, see Proposition \ref{prop:noc}.
Therefore, Proposition \ref{prop:t-hat-n} implies that the saddle points of first order can never be local minimizers. This observation explains why we use the term "saddle points" for them.
Among the saddle points of first order, we further distinguish those with vanishing and/or non-vanishing biactive multipliers.

\begin{definition}[Singular and regular saddle points of first order]
    Let $\bar x \in M$ be a saddle point of first order with multipliers $\lambda_{1,i}, \lambda_{2,i}$,$i \in I_{00}(\bar x)$. Then, it is called:
    \begin{itemize}
        \item singular if there exists $i \in I_{00}(\bar x)$ with either $\lambda_{1,i} =0$ or $\lambda_{2,i} =0$,
        \item regular if there exists $i \in I_{00}(\bar x)$ with both $\lambda_{1,i} \not =0$ and $\lambda_{2,i} \not =0$.
    \end{itemize} 
\end{definition}

Clearly, any saddle point of first order is singular, regular, or both. For the latter, note that their definitions depend on which pair of biactive multipliers is considered. It may well happen that with respect to one nonsmooth constraint a saddle point of first order is singular, wheres with respect to another it is regular. Additionally, singular saddle points of first order always violate ND2. The biactive index of regular saddle points of first order does not vanish, i.e. $BI \not =0$. Vice versa, if $BI \not =0$ at a saddle point of first order, then it is regular. 
These observations explain our choice of the terms "singular" and "regular" for the saddle points of first order. The notions of singular (in blue) and regular (in red) saddle points of first order are illustrated in Figure \ref{fig:spfo}.

We emphasize that a regular saddle point of first order admits two directions within the underlying nonsmooth structure, along which the objective function linearly decreases. A singular saddle points of first order can be made regular by arbitrarily small $C^1$-perturbations of the SNO defining functions. In order to illustrate these observations, we present examples of singular and regular saddle points of first order for MPCC. Similar examples can be analogously provided for other subclasses of SNO, such as MPVC, MPOC, MPSC, and MPDC.

\begin{example}[Regular saddle points of first order]
\label{ex:2}
     Consider the following MPCC:
\begin{equation}
  \label{eq:ex2}
  \min \,\,  \left(x_1-1\right)^2 + \left(x_2-1\right)^2 \quad \mbox{s.t.} \quad (x_1, x_2) \in \mathbb{K}_\text{c}.
\end{equation}
It is easy to see that $(0,0)$, $(1,0)$, and $(0,1)$ are the only T-stationary points for (\ref{eq:ex2}). More than that, both $(1,0)$ and $(0,1)$ solve (\ref{eq:ex2}). The biactive multipliers of $(0,0)$ are negative with $\lambda_{1}=-2$ and $\lambda_{2}=-2$. Hence, $(0,0)$ is a regular saddle point of first order. Actually, it is even nondegenerate with $BI=1$ and connects the minimizers $(1, 0)$ and $(0,1)$. Note that $(0,0)$ is not $N$-stationary for (\ref{eq:ex2}).
\end{example}

Qualitatively the same global structure as in Example \ref{ex:2} can be also induced by singular saddle points of first order.

\begin{example}[Singular saddle points of first order I]
\label{ex:1}
     Consider the following MPCC:
\begin{equation}
  \label{eq:ex1}
  \min \,\,  -x_1+\frac{1}{2}x_1^2 - x_2^2+\frac{1}{2}x_2^4 \quad \mbox{s.t.} \quad (x_1, x_2) \in \mathbb{K}_\text{c}.
\end{equation}
It is easy to see that $(0,0)$, $(1,0)$, and $(0,1)$ are the only T-stationary points for (\ref{eq:ex1}). More than that, both $(1,0)$ and $(0,1)$ solve (\ref{eq:ex1}). For the biactive multipliers of $(0,0)$ we have $\lambda_{1}=-1$ and $\lambda_{2}=0$. Hence, $(0,0)$ is a singular saddle point of first order. Although even degenerate, it connects the minimizers $(1,0)$ and $(0,1)$. Note that $(0,0)$ is $N$-, but not $\widehat N$-stationary for (\ref{eq:ex1}).
\end{example}

We point out that singular saddle points of first order do not need to lead to different local minimizers as it happens in Example \ref{ex:1}. Rather than that, they may lead just to one local minimizer. However, the case of regular saddle points of first order can be recovered here by applying sufficiently small, but arbitrary $C^1$-perturbations of the defining functions. 

\begin{example}[Singular saddle points of first order II]
\label{ex:4}
     Consider the following MPCC:
\begin{equation}
  \label{eq:ex4}
  \min \,\,  -x_1 +\frac{1}{2}x_1^2+ x_2^2-\frac{1}{2}x_2^4 \quad \mbox{s.t.} \quad (x_1, x_2) \in \mathbb{K}_\text{c}.
\end{equation}
It is easy to see that $(0,0)$, $(1,0)$, and $(0,1)$ are the only T-stationary points for (\ref{eq:ex4}). More than that, $(1,0)$ is the local minimizer of (\ref{eq:ex4}), whereas $(0,1)$ is its local maximizer. For the biactive multipliers of $(0,0)$ we have $\lambda_{1}=-1$ and $\lambda_{2}=0$. Hence, $(0,0)$ is a singular saddle point of first order. Note that  $(0,0)$ is again $N$-, but not $\widehat N$-stationary for (\ref{eq:ex4}). We see that here $(0,0)$ connects the minimizer $(1,0)$ with the maximizer $(0,1)$. 
Let us explain in which sense $(0,0)$ is a saddle point nevertheless. 
For that, let us consider the following perturbation of (\ref{eq:ex4}) with sufficiently small, but arbitrary $\varepsilon>0$:
\begin{equation}
  \label{eq:ex4a}
  \min \,\,  -x_1 +\frac{1}{2}x_1^2 - \varepsilon \cdot x_2+ x_2^2-\frac{1}{2}x_2^4 \quad \mbox{s.t.} \quad (x_1, x_2) \in \mathbb{K}_\text{c}.
\end{equation}
It is easy to see that the T-stationary points for (\ref{eq:ex4a}) are $(0,0)$, $(1,0)$, $(0,t_1)$, $(0,t_2)$, and $(0,t_3)$, where $0< t_1 < t_2 < t_3$ depend on $\varepsilon$. More than that, both $(1,0)$ and $(0,t_1)$ are local minimizers of (\ref{eq:ex4a}), $(0,t_2)$ is its turning point, and $(0,t_3)$ is its local maximizer. For the biactive multipliers of $(0,0)$ we have $\lambda_{1}=-1$ and $\lambda_{2}=-\varepsilon$. Hence, $(0,0)$ is a regular saddle point of first order. It connects two local minimizers $(1,0)$ and $(0,t_1)$ as in Example \ref{ex:2}.
\end{example}

It is not difficult to construct saddle points of higher order. Note that they are unstable with respect to $C^2$-perturbations of the defining functions, cf. \cite{jongen:2012}. Moreover, they are not relevant for our purpose of comparing different stationarity notions. 

\begin{example}[Saddle points of second order]
\label{ex:3}
     Consider the following MPCC:
\begin{equation}
  \label{eq:ex3}
  \min \,\,  -x_1^2+\frac{1}{2}x_1^4 - x_2^2+\frac{1}{2}x_2^4 \quad \mbox{s.t.} \quad (x_1, x_2) \in \mathbb{K}_\text{c}.
\end{equation}
It is easy to see that $(0,0)$, $(1,0)$, and $(0,1)$ are the only T-stationary points for (\ref{eq:ex3}). More than that, both $(1,0)$ and $(0,1)$ solve (\ref{eq:ex3}). For the biactive multipliers of $(0,0)$ we have $\lambda_{1}=0$ and $\lambda_{2}=0$. Hence, $(0,0)$ is not a saddle point of first order. Although even degenerate, it connects the minimizers $(1,0)$ and $(0,1)$. Note that $(0,0)$ is $\widehat N$-stationary for (\ref{eq:ex3}). 
\end{example}

It is also instructive to consider degenerate T-stationary points, which are not saddle points of first order. 

\begin{example}[Not saddle points of first order]
\label{ex:5}
     Consider the following MPCC:
\begin{equation}
  \label{eq:ex5}
  \min \,\,  x_1 -\frac{1}{2}x_1^2- x_2^2+\frac{1}{2}x_2^4 \quad \mbox{s.t.} \quad (x_1, x_2) \in \mathbb{K}_\text{c}.
\end{equation}
It is easy to see that $(0,0)$, $(1,0)$ and $(0,1)$ are the only T-stationary points for (\ref{eq:ex5}). More than that, $(0,1)$ is the local minimizer of (\ref{eq:ex5}), whereas $(1,0)$ is its local maximizer. For the biactive multipliers of $(0,0)$ we have $\lambda_{1}=1$ and $\lambda_{2}=0$. Hence, $(0,0)$ is not a saddle point of first order. Nevertheless, it connects the minimizer $(0,1)$ with the maximizer $(1,0)$ as in Example \ref{ex:4}, where we dealt with a singular saddle point of first order instead. The difference to the seemingly the same situation as in Example \ref{ex:4} can be again explained in terms of $C^1$-perturbations of the defining functions. Namely, by no means $(0,0)$ can be perturbed to a saddle point of first order. This is due to its positive multiplier $\lambda_1=1$. 
Let us additionally explain why $(0,0)$ is nevertheless $\widehat N$-stationary for (\ref{eq:ex5}). For that, let us consider the following MPCC:
\begin{equation}
  \label{eq:ex5a}
  \min \,\,  x_1 -\frac{1}{2}x_1^2 + x_2^2-\frac{1}{2}x_2^4 \quad \mbox{s.t.} \quad (x_1, x_2) \in \mathbb{K}_\text{c}.
\end{equation}
It is easy to see that $(0,0)$, $(1,0)$, and $(0,1)$ are the only T-stationary points for (\ref{eq:ex5a}). More than that, both $(1,0)$ and $(0,1)$ are maximizers of (\ref{eq:ex5a}), and $(0,0)$ is its local minimizer. For the biactive multipliers of $(0,0)$ we have here also $\lambda_{1}=1$ and $\lambda_{2}=0$. Although the multipliers of $(0,0)$ are the same as in case of (\ref{eq:ex5}), it connects here the maximizers $(0,1)$ and $(1,0)$.
We conclude that for both MPCCs (\ref{eq:ex5}) and (\ref{eq:ex5a}) we cannot distinguish the type of $(0,0)$ just on the signs of its biactive multipliers. The question on whether it connects minimizer and maximizer, on one side, or two maximizers, on the other side, cannot be resolved by considering first-order information of the objective functions. So, $(0,0)$ is $\widehat N$-stationary not only for (\ref{eq:ex5a}), which would be straightforward since it is a minimizer there, but also for (\ref{eq:ex5}).

\end{example}

Now, we are ready to adequately describe the differences between $N$ and $\widehat N$-stationarity, on one side, and between T- and $N$-stationarity, on the other side, cf.
Figures \ref{fig:stat1} and \ref{fig:spfo}.

\begin{proposition}[$N$- vs. $\widehat N$-stationarity]
  \label{prop:n-hat-n}
  An $N$-stationary point $\bar x\in M$ satisfying SNO-LICQ is not $\widehat N$-stationary for SNO if and only if it is a singular saddle point of first order.   
\end{proposition}

\begin{proposition}[T- vs. $N$-stationarity]
  \label{prop:t-n}
  A T-stationary point $\bar x\in M$ satisfying SNO-LICQ is not $N$-stationary for SNO if and only if it is a regular saddle point of first order.   
\end{proposition}

Finally, we turn our attention to the comparison between $\overline{N}$- and T-stationarity. In general, there are $\overline{N}$-, but not T-stationary points, see Figure \ref{fig:stat1}. We emphasize that they are irrelevant from the topologically point of view. Being stable with respect to $C^1$-perturbations of the defining functions, they cannot be transformed into local minimizers or saddle points.

\begin{example}[Non-T-stationary points]
\label{ex:6}
     Consider the following MPCC:
\begin{equation}
  \label{eq:ex6a}
  \min \,\,  x_1 -\frac{1}{2}x_1^2 - x_2+\frac{1}{2}x_2^2 \quad \mbox{s.t.} \quad (x_1, x_2) \in \mathbb{K}_\text{c}.
\end{equation}
It is easy to see that $(1,0)$ and $(0,1)$ are the only T-stationary points for (\ref{eq:ex6a}). More than that, $(0,1)$ is the local minimizer of (\ref{eq:ex6a}), whereas $(1,0)$ is its local maximizer. Note that $(0,0)$ is not T-stationary. Indeed, if (\ref{eq:grad}) holds for $(0,0)$, then only with the biactive multipliers $\lambda_{1}=1$ and $\lambda_{2}=-1$, a contradiction to T-stationarity. Nevertheless, $(0,0)$ can be viewed to connect the minimizer $(0,1)$ with the maximizer $(1,0)$ as in Examples \ref{ex:4} and \ref{ex:5}. Recall that in 
Example \ref{ex:4} we dealt with a singular saddle point of first order, i.e. being $N$-stationary, and in Example \ref{ex:5} the point under consideration was $\widehat N$-stationary. The difference to the seemingly the same situation as in Examples \ref{ex:4} and \ref{ex:5} can be once again explained in terms of $C^1$-perturbations of the defining functions. Namely, by no means $(0,0)$ can be perturbed to a saddle point of first order or to be confused with a local minimizer. This is due to the different signs of $\lambda_1=1$ and $\lambda_2=-1$. 
\end{example}

We express the observations from Propositions \ref{prop:t-hat-n}, \ref{prop:n-hat-n}, and \ref{prop:t-n} by means of the following formal equivalences:

\[
\begin{array}{rcl}
   \mbox{T-stationary} \,\backslash\,
   \widehat N\mbox{-stationary} & \Leftrightarrow &\mbox{saddle points of first order}, \\ \\
   N\mbox{-stationary} \,\backslash\,
   \widehat N\mbox{-stationary} & \Leftrightarrow & \mbox{singular saddle points of first order}, \\ \\   \mbox{T-stationary} \,\backslash\,
   N\mbox{-stationary} & \Leftrightarrow &\mbox{regular saddle points of first order}.
\end{array}
\]
Overall, the hierarchy of $\widehat N$-, $N$-, $\overline{N}$-, and $T$-stationarity notions can be interpreted in terms of saddle points. Starting by $\widehat N$-stationarity, one takes singular saddle points of first order in order to get $N$-stationarity. Afterwards, regular saddle points of first order enlarge $N$- to T-stationarity. It is worth to mention that nondegenerate saddle points with $BI \not =0$ are, in particular, regular saddle points of first order. They are not $N$-stationary, but play a crucial role within the Morse theory. Namely, if crossing the corresponding T-stationary level, the topology of the lower level set changes, see Section \ref{sec:t-stat}. As consequence, such nondegenerate saddle points with $BI \not =0$ may connect to different local minimizers, cf. Example \ref{ex:2}. Further, the enlargement of T- to $\overline{N}$-stationarity brings topologically irrelevant stationary points into consideration. We conclude that the advantage of T-stationarity compared to the geometrically motivated stationarity notions lies mainly in the incorporation of regular saddle points of first order. This allows to perform the topological study of SNO from the perspective of global optimization. 
$N$-stationarity, which adds only singular saddle points of first order to $\widehat N$-stationarity, is not enough, whereas $\overline{N}$-stationarity, which incorporates also non-T-stationary points, is too ambiguous for this purpose. $\widehat N$-stationarity admittedly contains all local minimizers, but misses the saddle points of first order, either singular or regular. 

\begin{figure}[h]
    \centering
    \includegraphics[trim=600 570 700 85, scale=0.75]{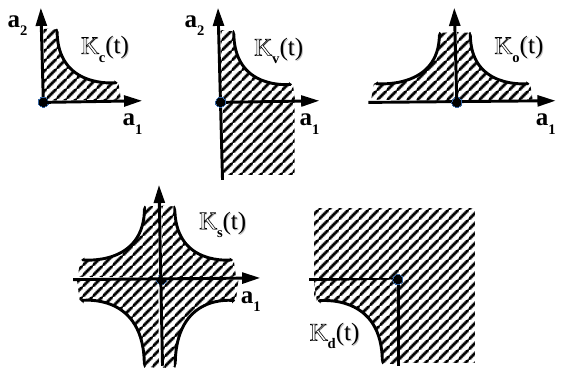}  
    \caption{Regularization of nonsmooth structures}
    \label{fig:smooth}
\end{figure}

The importance of T-stationarity from the practical point of view becomes clear if one regularizes SNO. The idea behind probably the most popular Scholtes type regularization is to associate SNO with nonlinear programming, see \cite{scholtes:2001}. This is done by substituting the cone $\mathbb{K}$ by $\mathbb{K}(t)$, where $t >0$ is a regularization parameter, see Figure \ref{fig:smooth}. Subsequently, Karush-Kuhn-Tucker points of the regularization are searched for: 
\[
\mbox{SNO}(t): \quad
\min_{x} \,\, f(x)\quad \mbox{s.\,t.} \quad x \in M(t)
\]
with
\[
    M(t)=\left\{x \in\R^n \left\vert\;
       \left(F_{1,i}(x), F_{2,i}(x)\right)\in \K(t), i=1,\ldots,m 
     \right.\right\}.
\]
It turns out that for $t \rightarrow 0$ the Karush-Kuhn-Tucker points of $\mbox{SNO}(t)$ converge to the T-stationary points of SNO. Moreover, the latter statement cannot be in general tightened by replacing the T-stationarity by any of the geometrically motivated notions of $\widehat N$- or $N$-stationarity. This is again due to the fact that they exclude regular saddle points of first order from consideration. Let us illustrate this phenomenon by means of a typical MPCC example, where for the regularization of $\K_\text{c}$ we set
\[
\K_\text{c}(t)=\left\{ (a_1,a_2) \in \R^2\, \left\vert\, a_1\cdot a_2 \leq t, a_1 \geq 0, a_2 \geq 0\right.\right\}.
\] 

\begin{example}[Scholtes type regularization]
\label{ex:2t}
     Consider the Scholtes type regularization of the MPCC from Example \ref{ex:2}:
\begin{equation}
  \label{eq:ex2t}
  \min \,\,  \left(x_1-1\right)^2 + \left(x_2-1\right)^2 \quad \mbox{s.t.} \quad (x_1, x_2) \in \mathbb{K}_\text{c}(t),
\end{equation}
where $t > 0$ is taken sufficiently small. It is easy to see that 
$(\sqrt{t},\sqrt{t})$, $\left(\frac{1+\sqrt{1-4t}}{2},\frac{1-\sqrt{1-4t}}{2}\right)$, and $\left(\frac{1-\sqrt{1-4t}}{2},\frac{1+\sqrt{1-4t}}{2}\right)$ are the only Karush-Kuhn-Tucker points for (\ref{eq:ex2t}). More than that, both $\left(\frac{1+\sqrt{1-4t}}{2},\frac{1-\sqrt{1-4t}}{2}\right)$ and $\left(\frac{1-\sqrt{1-4t}}{2},\frac{1+\sqrt{1-4t}}{2}\right)$ solve (\ref{eq:ex2t}). The saddle point $(\sqrt{t},\sqrt{t})$ with quadratic index one connects them. For $t \rightarrow 0$, the above Karush-Kuhn-Tucker points tend to the T-stationary points $(0,0)$, $(1,0)$, and $(0,1)$ of the original MPCC (\ref{eq:ex2}), respectively. Note that as the limit of $(\sqrt{t},\sqrt{t})$ we obtain $(0,0)$, which is neither $\widehat N$- nor $N$-stationary for (\ref{eq:ex2}). 
Recall that $(0,0)$ is a regular saddle point of first order for (\ref{eq:ex2}) instead. Thus, in order to fully describe the convergence properties of the regularization (\ref{eq:ex2t}), regular saddle points of first order cannot be neglected. 
\end{example}


\bibliography{sn-bibliography}

\end{document}